\documentclass[12pt,a4paper]{article}

\usepackage[latin1]{inputenc}
\usepackage[T1]{fontenc}
\usepackage[english]{babel}
\usepackage{amssymb,amsmath}

\newcommand{\augmentelargeur}[1]{
\addtolength{\evensidemargin}{-#1}
\addtolength{\oddsidemargin}{-#1}
\addtolength{\textwidth}{#1}
\addtolength{\textwidth}{#1}
}
\augmentelargeur{5mm}

\newcommand{\thmheadercommand}[1]{\textbf{\scshape{}#1}}

\def\d{{\mathrm{d}}}
\def\N{{\mathbb{N}}}
\def\Z{{\mathbb{Z}}}

\renewcommand{\geq}{\geqslant}
\renewcommand{\leq}{\leqslant}

\def\eps{\varepsilon}
\renewcommand{\epsilon}{\varepsilon}
\renewcommand{\phi}{\varphi}

\newcommand{\abs}[1]{\left|\mskip1mu#1\right|}

\newcommand{\presgroup}[2]{\left\langle\,#1 \mid  #2\,\right\rangle}

\newcounter{prop}
\newcounter{defi}
\newcounter{thm}
\newcounter{lem}

\newenvironment{dem}[1][]{\noindent{\thmheadercommand{Proof#1}}\,\,--\,\,}{$\square$\medskip}

\newenvironment{enonce}[1]{\medskip\noindent{\thmheadercommand{#1}}\,\,--\,\,\begin{slshape}}{\end{slshape}\medskip}

\newenvironment{enonce2}[1]{\medskip\noindent{\thmheadercommand{#1}}\,\,--\,\,}{\medskip}

\newenvironment{prop}[1][]{\refstepcounter{prop}
\begin{enonce}{Proposition \theprop{}#1}}{\end{enonce}}
\newenvironment{thm}[1][]{\refstepcounter{prop}
\begin{enonce}{Theorem \theprop{}#1}}{\end{enonce}}
\newenvironment{lem}[1][]{\refstepcounter{prop}
\begin{enonce}{Lemma \theprop{}#1}}{\end{enonce}}
\newenvironment{cor}[1][]{\refstepcounter{prop}
\begin{enonce}{Corollary \theprop{}#1}}{\end{enonce}}

\newenvironment{rem}[1][]{\refstepcounter{prop}
\begin{enonce2}{Remark \theprop{}#1}}{\end{enonce2}}

\newenvironment{defi*}[1][]{
\begin{enonce}{Definition#1}}{\end{enonce}}
\newenvironment{prop*}[1][]{
\begin{enonce}{Proposition#1}}{\end{enonce}}
\newenvironment{thm*}[1][]{
\begin{enonce}{Theorem#1}}{\end{enonce}}
\newenvironment{lem*}[1][]{
\begin{enonce}{Lemma#1}}{\end{enonce}}
\newenvironment{cor*}[1][]{
\begin{enonce}{Corollary#1}}{\end{enonce}}
\newenvironment{ex*}[1][]{
\begin{enonce}{Example#1}}{\end{enonce}}
\newenvironment{exo*}[1][]{
\begin{enonce2}{Exercise#1}}{\end{enonce2}}
\newenvironment{rem*}[1][]{
\begin{enonce2}{Remark#1}}{\end{enonce2}}

\usepackage{graphicx}
\usepackage[dvips]{color}

\title{Some small cancellation properties of random groups}

\author{Yann Ollivier}

\def\d{\partial}

\begin{document}

\maketitle

\begin{abstract}
We work in the density model of random groups.
We prove that they satisfy an isoperimetric inequality with sharp constant $1-2d$ depending
upon the density parameter $d$. This implies in particular a property generalizing
the ordinary $C'$ small cancellation condition, which could be termed
``macroscopic small cancellation''. This also sharpens the evaluation of
the hyperbolicity constant $\delta$.

As a consequence we get that the standard presentation of a random group
at density $d<1/5$ satisfies the Dehn algorithm and Greendlinger's Lemma,
and that it does not for $d>1/5$.
\end{abstract}

\section*{Statements}

Gromov introduced in~\cite{Gro93} the so-called density model of random
groups, which allows the study of generic groups with a very precise
control on the number of relators put in the group, depending on a
density parameter $d$.

A set of $m$ generators $a_1,\ldots,a_m$ being fixed, this model consists
in choosing a large length $\ell$ and a density parameter $0\leq d\leq
1$, and choosing at random a set $R$ of $(2m-1)^{d\ell}$ reduced words of
length $\ell$. The random group is then the group given by the
presentation $\presgroup{a_1,\ldots,a_m}{R}$. (Recall a word is \emph{reduced} if it
does not contain a generator immediately followed by its inverse).

In this model, we say that a property occurs \emph{with overwhelming
probability} if its probability of occurrence tends to $1$ as
$\ell\rightarrow\infty$ (everything else being fixed).

The basic intuition behind the model is that at density $d$, subwords of
length $(d-\eps)\ell$ of the relators will exhaust all possible reduced
words of this length. Also, at density $d$, with overwhelming probability
there are two relators sharing a subword of length $(2d-\eps)\ell$. We
refer to~\cite{Gro93}, \cite{Oll04} or~\cite{Ghy03} for a general
discussion on random groups and the density model.

The interest of this way to measure the number of relators in a
presentation is largely established by the following foundational theorem
of this theory, due to Gromov (\cite{Gro93}, see also~\cite{Oll04}).

\begin{thm}
If $d<1/2$, with overwhelming probability a random group at density $d$ is
infinite and hyperbolic.

If $d>1/2$, with overwhelming probability a random group at density $d$ is
either $\{e\}$ or $\Z/2\Z$.
\end{thm}

(Occurrence of $\Z/2\Z$ of course corresponds to even $\ell$.)

Other properties of random groups are known: Property~$T$ for $d>1/3$
\cite{Zuk03}, the small cancellation $C'(1/6)$ condition for $d<1/12$,
spectral properties of the random walk on the resulting group for $d<1/2$
\cite{Oll03rcg}, growth exponent for $d<1/2$ \cite{Oll03rg}, and several
properties at densities arbitrarily close to $0$ (see references
in~\cite{Ghy03} or~\cite{Oll04}).  The construction can be modified and
iterated in various ways to achieve specific goals \cite{Gro03}.

\bigskip

Hyperbolicity for $d<1/2$ is achieved by proving that van Kampen diagrams
satisfy some isoperimetric inequality (we refer to~\cite{LS77} for
definitions about van Kampen diagrams and to~\cite{Sho91} for the
equivalence between hyperbolicity and isoperimetry of van Kampen
diagrams). The main result of this paper is a sharp version of this
isoperimetric inequality.

\begin{thm}
\label{main}
For every $\eps>0$, with overwhelming probability,
every reduced van Kampen diagram $D$ in a random
group at density $d$ satisfies
\[
\abs{\d D}\geq (1-2d-\eps)\,\ell\abs{D}
\]
\end{thm}

This was already known to hold for diagrams of bounded size (see
Proposition~\ref{boundediso}), but the passage to all diagrams involves
the local-global hyperbolic principle of Gromov (see e.g.~\cite{Pap96}),
which implies a loss in the constants. After using this, the only
constant available for all diagrams was something like $(1-2d)/10^{20}$.

This inequality is sharp: indeed, at density $d$ there are very probably
two relators sharing a subword of length $(2d-\eps)\ell$, so that they
can be arranged to form a $2$-face van Kampen diagram of boundary length
$2(1-2d+\eps)\ell$. At density $d$ one can always glue some new relator
to any diagram along a path of length $(d-\eps)\ell$, so that adding
relators to this example provides an arbitrarily large diagram with the
same isoperimetric constant.

\bigskip

Besides its aesthetic interest as a sharp constant depending on density,
Theorem~\ref{main} also allows to prove several combinatorial properties
of the presentations obtained. First we give an immediate (but probably
unimportant) corollary having to do with small cancellation. Second, this
improves the estimate of the hyperbolicity constant $\delta$.
Theorem~\ref{main} also allows to prove that the critical density for
satisfaction of the Dehn algorithm and Greendlinger's Lemma is
$1/5$ (Theorem~\ref{dehn} below).
Finally, Theorem~\ref{main} will be used in~\cite{OW} to show that random
groups at densities $<1/6$ act on $CAT(0)$ cube complexes and satisfy the
Haagerup property.

Let us stress that the Dehn algorithm, as well as the properties studied
in~\cite{OW}, could not be obtained with the previous constant
$(1-2d)/10^{20}$, if only for the reason that this number is never
greater than $1/2$... So the improvement allows qualitative progress, not
only a quantitative one as for the hyperbolicity constant $\delta$.

\begin{cor}
\label{corsc}
For every $\eps>0$, with overwhelming probability, random groups at
density $d$ satisfy the following: Let $D_1$ and $D_2$ be two reduced van
Kampen diagrams and suppose that their boundaries share a common reduced
subword $w$. Suppose moreover that the diagram $D=D_1\cup_w D_2$ obtained
by gluing $D_1$ and $D_2$ along $w$ is reduced around $w$. Then we have
\[
\abs{w}\leq d\,(\abs{\d D_1}+\abs{\d D_2}) \,(1+\eps)
\]
\end{cor}

When $D_1$ and $D_2$ each consist of only one face, this exactly states
that random groups satisfy the $C'(2d)$ small cancellation property
(which implies hyperbolicity only when $d<1/12$). So this
property is a kind of ``macroscopic small cancellation''.

\begin{cor}
\label{cordelta}
At density $d$, for any $\eps>0$, with overwhelming probability the
hyperbolicity constant of a random group satisfies $\delta\leq 12\ell/(1-2d-\eps)^2$.
\end{cor}

Of course, this is not qualitatively different from the $10^{40}$ times
larger previous estimate.

\bigskip

Our last application of Theorem~\ref{main} has to do with the Dehn
algorithm and Greendlinger's Lemma, which are classical properties
considered in combinatorial group theory (see~\cite{LS77}). [More refs needed here!]

One might expect from Theorem~\ref{main} that the Dehn algorithm holds as
soon as $d<1/4$. Indeed, $d<1/4$ implies that some face of any reduced
diagram has at least $\ell/2$ boundary edges; but these might not be
consecutive. Actually the critical density is $1/5$.

\begin{thm}
\label{dehn}
If $d<1/5$, with overwhelming probability, the standard presentation of a
random group satisfies the Dehn algorithm and Greendlinger's Lemma.

More precisely, for any $\eps>0$, with overwhelming probability, in every
reduced van Kampen diagram with at least two faces, there are at least
two faces having more than $\frac{\ell}{2}+\frac{\ell}{2}(1-5d-\eps)$
consecutive edges on the boundary of the diagram.

If $d>1/5$, with overwhelming probability, the standard presentation of a
random group does not satisfy the Dehn algorithm nor Greendlinger's Lemma.
\end{thm}

This refers to the random presentation obtained by applying directly the
definition of the density model. In any $\delta$-hyperbolic group, the set
of words of length at most $8\delta$ representing the identity
constitutes a presentation of the group satisfying the Dehn
algorithm (\cite{Sho91}, Theorem~2.12); however, this set of words is
quite large, and computing it is feasible but tedious. Moreover it does
in general not satisfy the Greendlinger lemma.

%[Refs to Holt, Epstein etc here]

What happens at $d=1/5$ is not known (just as what happens for
infiniteness or triviality at $d=1/2$), but probably depends on more
precise subexponential terms in the number of relators of the
presentation, and so might not be very interesting.

\paragraph{Acknowledgements.} Part of the ideas presented here arose
during my stay in Montréal in July 2004 at the invitation of Daniel
T.~Wise, whom I would like to thank for helpful discussions and his so
warm welcome. I would also like to thank Thomas Delzant for having
insisted on the importance of the Dehn algorithm.

\section*{Proof of Theorem~\ref{main}}

We are going to prove Theorem~\ref{main} by bootstrapping on the
local-global principle. First, we recall the
result from~\cite{Gro93} (see also~\cite{Oll04}) on diagrams of bounded
size.

Suppose we are given a random presentation at density $d$, by reduced
relators of length $\ell$.

\begin{prop}
\label{boundediso}
For every $\eps>0$ and every $K\in \N$, with overwhelming probability,
every reduced van Kampen diagram with at most $K$ faces satisfies
\[
\abs{\d D}\geq (1-2d-\eps)\,\ell\abs{D}
\]
\end{prop}

Of course, the overwhelming probability is a priori not uniform in $K$
and $\eps$.

\begin{dem}
We only have to change a little bit the conclusion of the proof
in~\cite{Oll04}, p.~613. It is proven there that if $D$ is a reduced van
Kampen diagram involving $n\leq \abs{D}$ distinct relators
$r_1,\ldots,r_n$, with relator $r_i$ appearing $m_i$ times in the
diagram (we can assume $m_1\geq \ldots \geq m_n$), then there exist number $d_i$, $1\leq i\leq n$ such that:
\[
\abs{\d D}\geq (1-2d)\,\ell\abs{D}+2\sum d_i(m_i-m_{i+1})
\]
and such that the probability of this situation is at most $(2m)^{\inf
d_i}$ (\cite{Oll04}, p.~613). In particular, for fixed $\eps$, with overwhelming
probability we can suppose that $\inf d_i\geq -\ell\eps/2$.

If all $d_i$'s are non-negative, then we get $\abs{\d D}\geq
(1-2d)\,\ell\abs{D}$ as needed.

Otherwise, as $1\leq m_i\leq \abs{D}$ and $m_i\geq m_{i+1}$ we have $\sum d_i(m_i-m_{i+1})\geq
\abs{D}\inf d_i$ and so
\[
\abs{\d D}\geq (1-2d)\,\ell\abs{D}+2\abs{D}\inf d_i\geq
(1-2d-\eps)\,\ell\abs{D}
\]
\end{dem}

Then, using the local-global principle of hyperbolic geometry (see the
Proposition on page~613 of~\cite{Oll04}) we get that

\begin{prop}
\label{intermiso}
There exists a
constant $C>0$ such that any reduced van Kampen diagram $D$ (not
only those having at most $K$ faces) satisfy
\[
\abs{\d D}\geq C\,\ell\abs{D}
\]
\end{prop}

The constant $C$ is basically $1-2d$ divided by some huge constant (of
order $10^{20}$), so this is not what we need...

We solve the problem by a kind of bootstrapping: we will re-do some kind
of local-global principle to sharpen the constant, using the conclusions
of the above local-global principle.

\begin{prop}
\label{localglobaliso}
Suppose that for some $C>0$ any reduced van Kampen diagram $D$ satisfies \[
\abs{\d D}\geq C\,\ell\abs{D}
\]

Let $K\geq 3000/C^5$ and
suppose that any reduced van Kampen diagram $D$ with at most $K$ faces satisfies
\[
\abs{\d D}\geq \beta \,\ell\abs{D}
\]
Then any reduced van Kampen diagram satisfies
\[
\abs{\d D}\geq (\beta-14/\sqrt{KC})\,\ell\abs{D}
\]

\end{prop}

\begin{rem}
Here ``reduced'' could be replaced by ``having some property $P$'' with
$P$ a property such that any subdiagram of a diagram with $P$ also has
$P$. Indeed, in the proof we only use subdiagrams of a given diagram. In
other contexts (such as a version of this relative to a hyperbolic
initial group), this may be useful with $P$ ``being of minimal area'', or
``being stronly reduced with respect to a subpresentation'' (see
Definition~29 of~\cite{Oll04}, also compare the notion of graded
reducedness in~\cite{Ols91}).
\end{rem}

\begin{dem}

We need several lemmas.

\begin{lem}
\label{narrow}
Suppose that for some $C>0$ any reduced van Kampen diagram $D$ satisfies
\[
\abs{\d D}\geq C\,\ell\abs{D}
\]
(where we can suppose $C\leq 1$). Set $\alpha=1/\log(1/(1-C))\leq 1/C$.

Let $D$ be a reduced van Kampen diagram. Then
each face of $D$ is at distance at most $\alpha \log\abs{D}$ from the boundary of $D$.
\end{lem}

(A face adjacent to the boundary is said to be at distance $1$ from the
boundary, a face adjacent to such a face, at distance $2$, etc.)

\begin{dem}[ of the lemma]
We have $\abs{\d D}\geq C\ell\abs{D}$. So there are at least $C\abs{D}$
faces of $D$ adjacent to the boundary.

Applying the same reasoning to the (maybe not connected) diagram obtained
from $D$ by removing the boundary faces, we get by induction that the
number of faces of $D$ lying at distance at least $k$ from
the boundary is at most $(1-C)^{k-1}\abs{D}$. Taking $k=1+\alpha
\log\abs{D}$ (rounded up to the nearest integer) shows that there is less
than $1$ face at distance $k$ from the boundary.
\end{dem}

\begin{lem}
\label{cut}
Suppose that for some $C>0$ any reduced van Kampen diagram $D$ satisfies
\[
\abs{\d D}\geq C\,\ell\abs{D}
\]
(where we can suppose $C\leq 1$). Set $\alpha=1/\log(1/(1-C))\leq 1/C$.

Let $D$ be a reduced van Kampen diagram. Then $D$ can be partitioned into
two diagrams $D',D''$ by cutting it along a path of length at most
$2\alpha\ell\log\abs{D}$ such that each of $D'$ and $D''$ contains
at least one quarter of the boundary of $D$.
\end{lem}

\begin{dem}[ of the lemma]
By Lemma~\ref{narrow}, any face of $D$ lies at distance at most
$\alpha\log \abs{D}$ from the boundary.

Let $L$ be the boundary length of $D$ and mark four points $A,B,C,D$ on
$\d D$ at distance $L/4$ of each other. As $D$ is $\alpha
\log \abs{D}$-narrow, there exists a path of length at most
$2\alpha\ell\log \abs{D}$ joining either a point of $AB$ to a
point of $CD$
or a point of $AD$ to a point of $BC$, which provides the desired
cutting.
\end{dem}

This allows to prove one step of the local-global passage.

\begin{lem}
\label{steplocalglobal}
Suppose that for some $C>0$ any reduced van Kampen diagram $D$ satisfies \[
\abs{\d D}\geq C\,\ell\abs{D}
\]
and that, for some $A\geq 3000/C^4$, any reduced van Kampen diagram $D$ with boundary length at most $A\ell$ satisfies
\[
\abs{\d D}\geq \beta \,\ell\abs{D}
\]

Then any reduced van Kampen diagram $D$ with boundary length at most $7A\ell/6$ satisfies
\[
\abs{\d D}\geq (\beta-1/\sqrt{A})\,\ell\abs{D}
\]
\end{lem}

\begin{dem}[ of the lemma]
Let $D$ be a reduced van Kampen diagram of boundary length between
$A\ell$ and $7A\ell/6$. We have $A\leq \abs{D}\leq 7A/6C$.

By Lemma~\ref{cut}, we can partition $D$ into two diagrams $D'$ and
$D''$, each of them containing at least one quarter of the boundary
length of $D$. So we have $\abs{\d D'}\leq 3\abs{\d
D}/4+2\ell\alpha\log \abs{D}\leq \ell(7A/8+2\alpha\log(7A/6C))$ and
likewise for $D''$.

Choose $A$ large enough (depending only $C$) so that
$2\alpha\log(7A/6C)\leq A/8$. Then both $D'$ and $D''$ have boundary
length at most $A\ell$. So we have
\[
\abs{\d D'}\geq \beta\ell\abs{D'}\ \text{and}\ \abs{\d D''}\geq
\beta\ell\abs{D''}
\]

Now we choose $A$ large enough (depending on $C$) so that
$4\alpha\log(7A/6C))\leq\sqrt{A}$ (taking $A=3000/C^4$ is enough). We have
\begin{align*}
\abs{\d D}&=\abs{\d D'}+\abs{\d D''}-2\abs{\d D' \cap \d D''}
\\&\geq
\abs{\d D'}+\abs{\d D''}-4\ell\alpha\log\abs{D}
\\&\geq
\beta\,\ell(\abs{D'}+\abs{D''})-\ell\sqrt{A}
\\&\geq (\beta-1/\sqrt{A})\,\ell\abs{D}
\end{align*}
since $\abs{D}\geq A$.
\end{dem}

Now we are ready to prove Proposition~\ref{localglobaliso}. Take $K=A/C$
where $A\geq 3000/C^4$. Then any reduced van
Kampen diagram of boundary length at most $A\ell$ has at most $K$ faces,
so that the assumption of Proposition~\ref{localglobaliso} implies the
assumption of Lemma~\ref{steplocalglobal}.

So applying this lemma,
we know that when going from diagrams of size $A\ell$ to diagrams of size
$7A\ell/6$, the isoperimetric constants $\beta$ worsens by
$-1/\sqrt{A}$. So by induction we are able to show that diagrams $D$ of size
between $A\ell(7/6)^k$ and $A\ell(7/6)^{k+1}$ satisfy the isoperimetric
inequality $\abs{\d D}\geq \beta_k \,\ell\abs{D}$ with
$\beta_0=\beta-1/\sqrt{A}$ and
$\beta_{k+1}=\beta_k-\frac{1}{\sqrt{A(7/6)^k}}$
so that $\beta_k\geq \beta-14/\sqrt{A}$ for any $k$.
This proves Proposition~\ref{localglobaliso}.
\end{dem}

Now the proof of Theorem~\ref{main} is clear: take the isoperimetric
constant $C$ provided by Proposition~\ref{intermiso}. Take $K$ so that
$K\geq 3000/C^5$ and $14/\sqrt{KC}\leq \eps/2$. By
Proposition~\ref{boundediso}, with overwhelming probability, we can
suppose that any reduced van Kampen $D$ diagram with at most $K$ faces
satisfies $\abs{\d D}\geq (1-2d-\eps/2)\,\ell\abs{D}$. Now apply
Proposition~\ref{localglobaliso} to end the proof.

Corollary~\ref{corsc} is easy. Let $D=D_1\cup_w D_2$. 
Since $\abs{\d D}\geq (1-2d-\eps)\,\ell\abs{D}$, the number of internal edges
of $D$ is at most $(d+\eps/2)\ell\abs{D}$. So a fortiori $\abs{w}\leq (d+\eps/2)\ell\abs{D}$.
Now
\begin{align*}
\abs{w}&\leq (d+\eps/2)\,\ell\abs{D}\leq\frac{d+\eps}{1-2d-\eps}\,\abs{\d D}
\\&= \frac{d+\eps/2}{1-2d-\eps}\left(\abs{\d D_1}+\abs{\d D_2}-2\abs{w}\right)
\end{align*}
and so
\[
\abs{w}\leq (d+\eps/2)\,\left(\abs{\d D_1}+\abs{\d D_2}\right)
\]
as needed.

Corollary~\ref{cordelta} is obtained by applying Proposition~7
of~\cite{Oll03rg} (which is only Theorem~2.5 of~\cite{Sho91} where we
took care of the constants).

\section*{The Dehn algorithm and Greendlinger's Lemma}

We now turn to the proof of Theorem~\ref{dehn}. Since the Greendlinger
Lemma is stronger than the Dehn algorithm, it suffices to prove the
former for $d<1/5$ and disprove the latter for $d>1/5$.

\paragraph{Greendlinger's Lemma for $d<1/5$.}
We begin by a lemma which is weaker in the sense that we do not ask for
the boundary edges to be consecutive. We will then conclude by a standard
argument.

\begin{lem}
\label{badgreendlinger}
For any $\eps>0$, with overwhelming probability, at density $d$ the
following holds:

Let $D$ be a reduced van Kampen diagram with at least two faces.
There exist two faces of $D$ each having at least
$\ell(1-5d/2-\eps)$ edges on the boundary of $D$ (maybe not consecutive).

\end{lem}

Observe that when $d<1/5$ this is more than $\ell/2$ (for small enough
$\eps$ depending on $1/5-d$).
This lemma is also valid at densities larger than
$1/5$ but becomes trivial at $d=2/5$.

\begin{dem}[ of the lemma]
Let $D$ be a reduced van Kampen diagram with at least two faces.

Let $f$ be a face of $D$ having the greatest number of edges on the
boundary. Say $f$ has $\alpha\ell$ edges on the boundary. Suppose that
any face other than $f$ has no more than $\beta\ell$ edges on the
boundary. We want to show that $\beta\geq 1-5d/2-\eps$. So suppose that
$\beta < 1-5d/2-\eps$. (The reader may find more convenient to read the
following skipping the $\eps$'s.)

Consider also the (maybe not connected, but this does not matter) diagram
$D'$ obtained by removing face $f$ from $D$. We have $\abs{\d D'}=\abs{\d
D}+\ell-2\alpha\ell$.

By definition of $\alpha$ and $\beta$ we have $\abs{\d D}\leq
\beta\ell(\abs{D}-1)+\alpha\ell$ and consequently $\abs{\d D'}\leq
\beta\ell(\abs{D}-1)+\ell-\alpha\ell$.

But by Theorem~\ref{main}, with overwhelming probability we can suppose
that we have $\abs{\d D}\geq
(1-2d-\eps/2)\,\ell\abs{D}$ and $\abs{\d D'}\geq
(1-2d-\eps/2)\,\ell\abs{D'}=(1-2d-\eps/2)\,\ell\,(\abs{D}-1)$. So combining
these inequalities we get
\begin{align*}
(1-2d-\eps/2)\abs{D} &\leq \beta(\abs{D}-1)+\alpha
\\
(1-2d-\eps/2)\,(\abs{D}-1) &\leq \beta(\abs{D}-1)+1-\alpha
\end{align*}
or, since we assumed that $\beta<1-5d/2-\eps$,
\begin{align*}
(1-2d-\eps/2)\abs{D} & < (1-5d/2-\eps)\,(\abs{D}-1)+\alpha
\\
(1-2d-\eps/2)\,(\abs{D}-1) & < (1-5d/2-\eps)\,(\abs{D}-1)+1-\alpha
\end{align*}
which yield respectively
\begin{align}
\label{isopDbeta}
\abs{D} &< \frac{\alpha+5d/2-1+\eps}{d/2+\eps/2}
\\
\label{isopD'beta}
\abs{D} &<\frac{d/2+1-\alpha+\eps/2}{d/2+\eps/2}
\end{align}

Either $\alpha\leq 1-d-\eps/4$ or $\alpha\geq 1-d-\eps/4$. In any case,
one of (\ref{isopDbeta}) or (\ref{isopD'beta}) gives
\[
\abs{D} < \frac{3d/2+3\eps/4}{d/2+\eps/2}<3
\]
(generally, a face having more than $(1-d)\ell$ on the boundary is the
frontier at which it is more interesting to remove this face before
applying Theorem~\ref{main}).

The case $\abs{D}\leq 2$ is easily ruled out. So we get a contradiction,
and the lemma is proven.
\end{dem}

This somewhat obscure proof and the role of $1/5$ will become clearer in
the next paragraph, when we will build a $3$-face diagram for
$d>1/5$ with only one face having more than $\ell/2$ boundary edges.

% 
% 
% We now turn to the proof of Corollary~\ref{cordehn}. Let $d<1/6$. By
% Theorem~\ref{main}, setting $\eps=1/6-d$, we can assume that
% $\abs{D}\geq(1-2d-\eps)\,\ell\abs{D}=(2/3+\eps)\,\ell\abs{D}$ for any
% reduced van Kampen diagram $D$.
% 
% Let $D$ be a reduced van Kampen diagram with at least two faces. We want
% to prove that there are two faces of $D$ having more than $\ell/2$
% consecutive edges on the boundary of $D$.
% 
% It is clear that there is some face having more than $\ell/2$ edges on
% the boundary of $D$ (but maybe not consecutive). Indeed, otherwise we
% would have $\abs{\d D}\leq \abs{D}\ell/2$, in contradiction with
% $\abs{\d D}\geq (2/3+\eps)\,\ell\abs{D}$.
% 
% We can even prove that there are two such faces. Indeed, if any face but
% one has no more than $\ell/2$ edges on the boundary we have
% $\abs{\d D}\leq (\abs{D}-1)\ell/2+\ell$, which, combined with $\abs{\d
% D}\geq (2/3+\eps)\,\ell\abs{D}$, yields $\abs{D}\leq 3/(1+\eps)$ so that
% $\abs{D}\leq 2$, an easily ruled out case.
% 

Back to the proof of Greendlinger's Lemma for $d<1/5$.
If we face a diagram $D$ such that the intersection of the boundary of
any face of $D$ with the boundary of $D$ is connected, then
Lemma~\ref{badgreendlinger} provides what we want.

Now we apply a standard argument to prove that this case is enough.
Suppose that some face of $D$ has a non-connected intersection with
the boundary, having two (or more) boundary components, so that this
face separates the rest of the diagram into two (or more) components.
Call \emph{good} a face having exactly one boundary component and
\emph{bad} a face with two or more boundary components (there are also
internal faces, which we are not interested in).

Decompose $D$ into bad faces and maximal parts without bad faces. Call
such a maximal part \emph{extremal} if it is in contact with only one bad
face. It is clear that there are at least two such extremal parts.

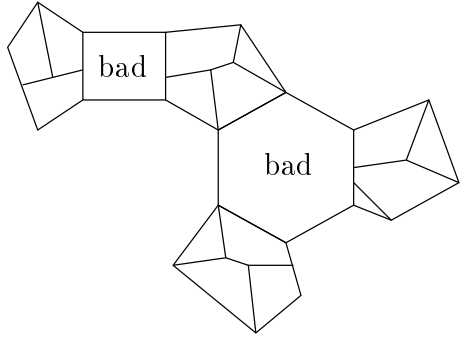
\begin{center}
\begin{picture}(0,0)%
\includegraphics{badfaces.pstex}%
\end{picture}%
\setlength{\unitlength}{4144sp}%
\begingroup\makeatletter\ifx\SetFigFont\undefined
% extract first six characters in \fmtname
\def\x#1#2#3#4#5#6#7\relax{\def\x{#1#2#3#4#5#6}}%
\expandafter\x\fmtname xxxxxx\relax \def\y{splain}%
\ifx\x\y   % LaTeX or SliTeX?
\gdef\SetFigFont#1#2#3{%
  \ifnum #1<17\tiny\else \ifnum #1<20\small\else
  \ifnum #1<24\normalsize\else \ifnum #1<29\large\else
  \ifnum #1<34\Large\else \ifnum #1<41\LARGE\else
     \huge\fi\fi\fi\fi\fi\fi
  \csname #3\endcsname}%
\else
\gdef\SetFigFont#1#2#3{\begingroup
  \count@#1\relax \ifnum 25<\count@\count@25\fi
  \def\x{\endgroup\@setsize\SetFigFont{#2pt}}%
  \expandafter\x
    \csname \romannumeral\the\count@ pt\expandafter\endcsname
    \csname @\romannumeral\the\count@ pt\endcsname
  \csname #3\endcsname}%
\fi
\fi\endgroup
\begin{picture}(2724,2004)(394,-1963)
\put(946,-421){\makebox(0,0)[lb]{\smash{\SetFigFont{12}{14.4}{rm}{\color[rgb]{0,0,0}bad}%
}}}
\put(1936,-1006){\makebox(0,0)[lb]{\smash{\SetFigFont{12}{14.4}{rm}{\color[rgb]{0,0,0}bad}%
}}}
\end{picture}

\end{center}

To reach the conclusion it is enough to find in any extremal part a good
face having more than $\ell(1-5d/2-\eps)$ edges on the boundary.  So let
$f$ be a bad face in contact with an extremal part $P$ with no bad faces.

Consider the diagram $D'=P\cup f$. This diagram has no bad faces now, and
so there are two faces in it having more than $\ell(1-5d/2-\eps)$
consecutive edges on the boundary. One of these may be $f$, but the other
one has to be in $P$ and so has more than $\ell(1-5d/2-\eps)$ consecutive
edges on the boundary of $D$ as well.

\paragraph{A counter-example for $d>1/5$.} Here we
show that the presentation does not satisfy the Dehn algorithm as soon as
$d>1/5$. So take $d>1/5$ and fix some $\eps>0$.

We can with overwhelming probability find two relators
$r_1,r_2$ sharing a common subword $w$ of length $(2d-\eps)\ell$. Once those are
chosen, let $x$ be the subword of length $(d-\eps)\ell$ of the boundary
of the diagram
$r_1\cup_w r_2$ occurring around the $w$-gluing and having length
$(d-\eps)\ell/2$ on each side of the $w$-gluing (see picture below). (When $d>2/5$ there is less
than this left on the boundary of $r_1\cup_w r_2$; but the situation is even
easier at larger densities and so we leave this detail aside).

At density $d$, subwords of length $(d-\eps)\ell$ of the relators exhaust
all reduced words of length $(d-\eps)\ell$. So it is possible to find a relator
$r_3$ gluing to $r_1\cup_w r_2$ along $x$. After this
operation $r_1$ and $r_2$ each have less than
$1-(2d-\eps)\ell-(d/2-\eps/2)\ell=(1-5d/2+3\eps/2)\ell$ of their length on the
boundary (see the picture below), which is less than $\ell/2$ when
$d>1/5$, for small enough $\eps$. Compare Lemma~\ref{badgreendlinger} ---
which is thus sharp.

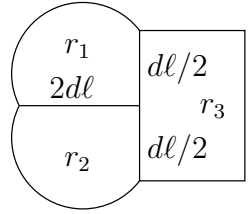
\begin{center}
\begin{picture}(0,0)%
\includegraphics{notdehn1_5.pstex}%
\end{picture}%
\setlength{\unitlength}{4144sp}%
\begingroup\makeatletter\ifx\SetFigFont\undefined
% extract first six characters in \fmtname
\def\x#1#2#3#4#5#6#7\relax{\def\x{#1#2#3#4#5#6}}%
\expandafter\x\fmtname xxxxxx\relax \def\y{splain}%
\ifx\x\y   % LaTeX or SliTeX?
\gdef\SetFigFont#1#2#3{%
  \ifnum #1<17\tiny\else \ifnum #1<20\small\else
  \ifnum #1<24\normalsize\else \ifnum #1<29\large\else
  \ifnum #1<34\Large\else \ifnum #1<41\LARGE\else
     \huge\fi\fi\fi\fi\fi\fi
  \csname #3\endcsname}%
\else
\gdef\SetFigFont#1#2#3{\begingroup
  \count@#1\relax \ifnum 25<\count@\count@25\fi
  \def\x{\endgroup\@setsize\SetFigFont{#2pt}}%
  \expandafter\x
    \csname \romannumeral\the\count@ pt\expandafter\endcsname
    \csname @\romannumeral\the\count@ pt\endcsname
  \csname #3\endcsname}%
\fi
\fi\endgroup
\begin{picture}(1417,1255)(576,-1138)
\put(811,-466){\makebox(0,0)[lb]{\smash{\SetFigFont{12}{14.4}{rm}{\color[rgb]{0,0,0}$2d\ell$}%
}}}
\put(901,-196){\makebox(0,0)[lb]{\smash{\SetFigFont{12}{14.4}{rm}{\color[rgb]{0,0,0}$r_1$}%
}}}
\put(901,-871){\makebox(0,0)[lb]{\smash{\SetFigFont{12}{14.4}{rm}{\color[rgb]{0,0,0}$r_2$}%
}}}
\put(1396,-331){\makebox(0,0)[lb]{\smash{\SetFigFont{12}{14.4}{rm}{\color[rgb]{0,0,0}$d\ell/2$}%
}}}
\put(1396,-826){\makebox(0,0)[lb]{\smash{\SetFigFont{12}{14.4}{rm}{\color[rgb]{0,0,0}$d\ell/2$}%
}}}
\put(1711,-556){\makebox(0,0)[lb]{\smash{\SetFigFont{12}{14.4}{rm}{\color[rgb]{0,0,0}$r_3$}%
}}}
\end{picture}

\end{center}

Note for later use that at this step, the boundary length of the diagram
so obtained is $(3-6d+4\eps)\ell$. This is the smallest possible value
compatible with Theorem~\ref{main}, up to the $\eps$'s.

But (thanks to the $\eps$'s) this will not only happen once but
arbitrarily many times as $\ell\rightarrow\infty$, so we can find another
independent triple $(r'_1,r'_2,r'_3)$ giving rise to the same
configuration.

Now if $r_3$ and $r'_3$ share only a single letter in the region of
length $\ell/5$ opposite to the position where they glue to $r_1\cup_w
r_2$ (resp.\ $r'_1\cup_{w'} r'_2$), and this happens all the time, then
we can form a diagram in which $r_3$ and $r'_3$ become faces having no
more than $\ell/2$ consecutive edges on the boundary (they are bad faces
in the terminology of the previous proof). So if $d>1/5$, no face of this
diagram has more than $\ell/2$ consecutive edges on the boundary
(although the two bad faces have more than $\ell/2$ non-consecutive
boundary edges).

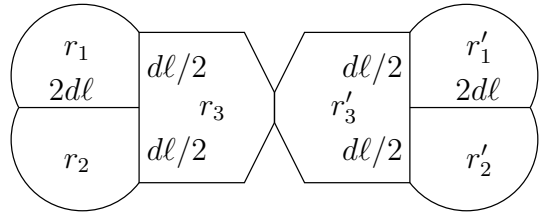
\begin{center}
\begin{picture}(0,0)%
\includegraphics{notdehn1_5full.pstex}%
\end{picture}%
\setlength{\unitlength}{4144sp}%
\begingroup\makeatletter\ifx\SetFigFont\undefined
% extract first six characters in \fmtname
\def\x#1#2#3#4#5#6#7\relax{\def\x{#1#2#3#4#5#6}}%
\expandafter\x\fmtname xxxxxx\relax \def\y{splain}%
\ifx\x\y   % LaTeX or SliTeX?
\gdef\SetFigFont#1#2#3{%
  \ifnum #1<17\tiny\else \ifnum #1<20\small\else
  \ifnum #1<24\normalsize\else \ifnum #1<29\large\else
  \ifnum #1<34\Large\else \ifnum #1<41\LARGE\else
     \huge\fi\fi\fi\fi\fi\fi
  \csname #3\endcsname}%
\else
\gdef\SetFigFont#1#2#3{\begingroup
  \count@#1\relax \ifnum 25<\count@\count@25\fi
  \def\x{\endgroup\@setsize\SetFigFont{#2pt}}%
  \expandafter\x
    \csname \romannumeral\the\count@ pt\expandafter\endcsname
    \csname @\romannumeral\the\count@ pt\endcsname
  \csname #3\endcsname}%
\fi
\fi\endgroup
\begin{picture}(3169,1255)(576,-1138)
\put(811,-466){\makebox(0,0)[lb]{\smash{\SetFigFont{12}{14.4}{rm}{\color[rgb]{0,0,0}$2d\ell$}%
}}}
\put(901,-196){\makebox(0,0)[lb]{\smash{\SetFigFont{12}{14.4}{rm}{\color[rgb]{0,0,0}$r_1$}%
}}}
\put(901,-871){\makebox(0,0)[lb]{\smash{\SetFigFont{12}{14.4}{rm}{\color[rgb]{0,0,0}$r_2$}%
}}}
\put(1396,-331){\makebox(0,0)[lb]{\smash{\SetFigFont{12}{14.4}{rm}{\color[rgb]{0,0,0}$d\ell/2$}%
}}}
\put(1396,-826){\makebox(0,0)[lb]{\smash{\SetFigFont{12}{14.4}{rm}{\color[rgb]{0,0,0}$d\ell/2$}%
}}}
\put(1711,-556){\makebox(0,0)[lb]{\smash{\SetFigFont{12}{14.4}{rm}{\color[rgb]{0,0,0}$r_3$}%
}}}
\put(2656,-556){\makebox(0,0)[rb]{\smash{\SetFigFont{12}{14.4}{rm}{\color[rgb]{0,0,0}$r'_3$}%
}}}
\put(2926,-331){\makebox(0,0)[rb]{\smash{\SetFigFont{12}{14.4}{rm}{\color[rgb]{0,0,0}$d\ell/2$}%
}}}
\put(2926,-826){\makebox(0,0)[rb]{\smash{\SetFigFont{12}{14.4}{rm}{\color[rgb]{0,0,0}$d\ell/2$}%
}}}
\put(3466,-196){\makebox(0,0)[rb]{\smash{\SetFigFont{12}{14.4}{rm}{\color[rgb]{0,0,0}$r'_1$}%
}}}
\put(3466,-871){\makebox(0,0)[rb]{\smash{\SetFigFont{12}{14.4}{rm}{\color[rgb]{0,0,0}$r'_2$}%
}}}
\put(3511,-466){\makebox(0,0)[rb]{\smash{\SetFigFont{12}{14.4}{rm}{\color[rgb]{0,0,0}$2d\ell$}%
}}}
\end{picture}

\end{center}

This is not enough to disprove the Dehn algorithm: this algorithm only
demands that for any reduced word representing $e$, there exists some van
Kampen diagram with the boundary face property. There could exist another
van Kampen diagram with the same boundary word as $D$, in which some face
would have more than $\ell/2$ consecutive edges on the boundary. So let
$r_4$ be this face; this means that we can glue $r_4^{-1}$ to the
previous diagram $D$ to get a new diagram $D'$ with $7$ faces; since
$r_4$ has more than half of its length on the boundary we have $\abs{\d
D'}<\abs{\d D}$.

Either $D'$ is reduced or $r_4$ is equal to some relator $r_i$ already present
in the diagram.

In the latter case, this means that we can glue a copy of $r_i^{-1}$
along $r_i$ on the boundary of the diagram $D$ along more than $\ell/2$
edges. But this means that before gluing $r_i^{-1}$ we could have folded
some letters of $r_i$ with neighbouring letters in the boundary of $D$.
This is excluded if we assume (as we can always do) that the boundary of
$D$ is reduced.

In the former case when $D'$ is reduced, using what we noted above we get
that $\abs{\d D}=(3-6d+4\eps)\ell\times2-2= 6(1-2d)\,\ell+8\eps\ell-2$.
Since $\abs{\d D'} < \abs{\d D}$ we get $\abs{\d D'}<
6(1-2d)\,\ell+8\eps\ell-2$. But by Theorem~\ref{main}, for any $\eps'$ we
have $\abs{\d D'}\geq 7(1-2d-\eps')\ell$, which is a contradiction for
small enough values of $\eps$ and $\eps'$.

\end{document}